\documentclass{amsproc}
\usepackage{color,graphicx,times}

\usepackage{amssymb}
\usepackage{amsmath}
\usepackage{amscd}
\usepackage{float}
\usepackage{graphicx}
\usepackage{amsfonts}
\usepackage{pb-diagram}
\usepackage{eufrak}

\theoremstyle{definition}

\theoremstyle{remark}

\numberwithin{equation}{section}

\begin{document}

\newcommand{\CrossStarGlyph}{\raisebox{-0.25\height}{\includegraphics[width=0.5cm]{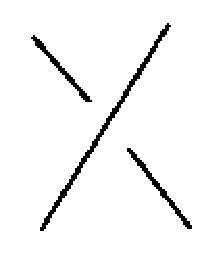}}}
\newcommand{\YMGlyph}{\raisebox{-0.25\height}{\includegraphics[width=0.5cm]{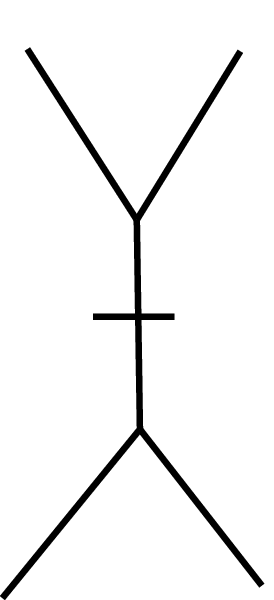}}}

\newcommand{\YMTens}{\raisebox{-0.25\height}{\includegraphics[width=0.5cm]{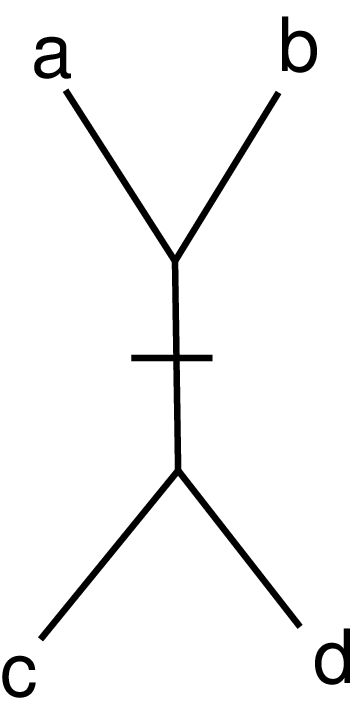}}}

\newcommand{\CDiag}{\raisebox{-0.25\height}{\includegraphics[width=0.5cm]{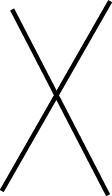}}}
\newcommand{\VDiag}{\raisebox{-0.25\height}{\includegraphics[width=0.5cm]{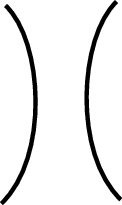}}}
\newcommand{\CDotDiag}{\raisebox{-0.25\height}{\includegraphics[width=0.5cm]{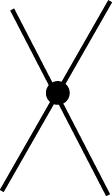}}}
\newcommand{\CDotDa}{\raisebox{-0.25\height}{\includegraphics[width=0.5cm]{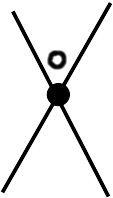}}}
\newcommand{\CDotDb}{\raisebox{-0.25\height}{\includegraphics[width=0.5cm]{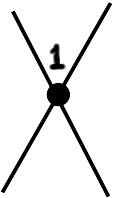}}}

\title{ Categorification of Chromatic, Dichromatic and Penrose Polynomials}

\author{Louis H. Kauffman}
\address{Department of Mathematics, Statistics and Computer Science, University of Illinois at Chicago, 851 South Morgan Street, Chicago Illinois, 60607-7045 $and$
International Institute for Sustainability with Knotted Chiral Meta Matter (WPI-SKCM2),
Hiroshima University, 1-3-1 Kagamiyama, Higashi-Hiroshima, Hiroshima 739-8526, Japan}
\email{loukau@gmail.com}




\begin{abstract}
This paper discusses ways to categorify chromatic, dichromatic and Penrose polynomials, including categorifications of integer evaluations of chromatic polynomials. 
 We show that with an appropriate choice of variables the coefficients of the Potts partition function at different energy levels are given by Euler characteristics of appropriate parts of a bigraded homology theory associated with the model. In the case of the dichromatic polynomial for graphs, we show that the two variable polynomial can be seen as a sum of powers of one variable multiplied by coefficients that are ``impropriety" coloring polynomials for the underlying graph. An impropriety polynomial $C_{G}^{i}(n)$ counts the number of colorings in $n$ colors  of the graph that are not proper at a given number $i$  of  edges in the graph. 
The last section of the paper categorifies coloring evaluations rather than coloring polynomials. We then obtain a range of possible chain complexes and homology theories such that the chromatic evaluation is equal to the Euler characteristic of the homology. The freedom of choice in making such chain complexes is related to possible associative algebra structures on the set of colors.
\end{abstract}

\maketitle

\section{Introduction}
In this paper we discuss ways to categorify chromatic, dichromatic and Penrose polynomials, including integer evaluations of chromatic polynomials. 
Constructions are motivated by examining the shape of state expansions of the corresponding polynomials with an eye 
toward making them into graded Euler characteristics of an associated chain complex and its homology theory.\\

The purpose of this paper is to give the definitions and constructions for these homology theories. A subsequent paper will consider calculations and relationships with coloring problems for graphs. In the case of the classical dichromatic polynomial for a graph, there is a relationship with the Potts model in statistical mechanics due to Neville Temperley (See \cite{Baxter,Temperley,K,LKStat,KP,DKT}). We show that with an appropriate choice of variables the coefficients of the Potts partition function at different energy levels are given by Euler characteristics of appropriate parts of a bigraded homology theory associated with the model. Also in the case of the dichromatic polynomial for graphs, we show that the two variable polynomial can be seen as a sum of powers of one variable multiplied by coefficients that are ``impropriety" coloring polynomials for the underlying graph. An impropriety polynomial $C_{G}^{i}(n)$ counts the number of colorings in $n$ colors  of the graph that are not proper at a given number $i$  of  edges in the graph. At these improper edges the coloring gives the same color to the ends of the edge. When $i=0$ then $C_{G}^{0}(n) = C_{G}(n)$ counts the number of proper colorings of the graph where distinct colors occur at the ends of every edge in the graph. The last section of the paper categorifies coloring evaluations rather than coloring polynomials. We then obtain a range of possible chain complexes and homology theories such that the chromatic evaluation is equal to the Euler characteristic of the homology (without extra grading). The freedom of choice in making such chain complexes is related to associative algebra structures on the set of colors.\\

In Section 2 the paper discusses categorification of the classical chromatic polynomial, starting with the Whitney logical identity for coloring as a deletion-contraction formula. It follows from the Whitney formula that the 
number of proper vertex colorings (vertices connected by an edge properly have distinct colors)  of a graph $G$, using $n$ colors, $C_{G}(n),$  is given  by the formula
$$C_{G}(n) = \sum_{S} (-1)^{e(S)}n^{c(S)}$$ where S runs over subgraphs of $G$ with $e(s)$ the number of edges of $S$ and $c(S)$ the number of connected components of $S.$  In this summation each subgraph consists in all the nodes of $G$ and a subset of the edges of $G.$\\

\begin{figure}
\includegraphics[width=8cm]{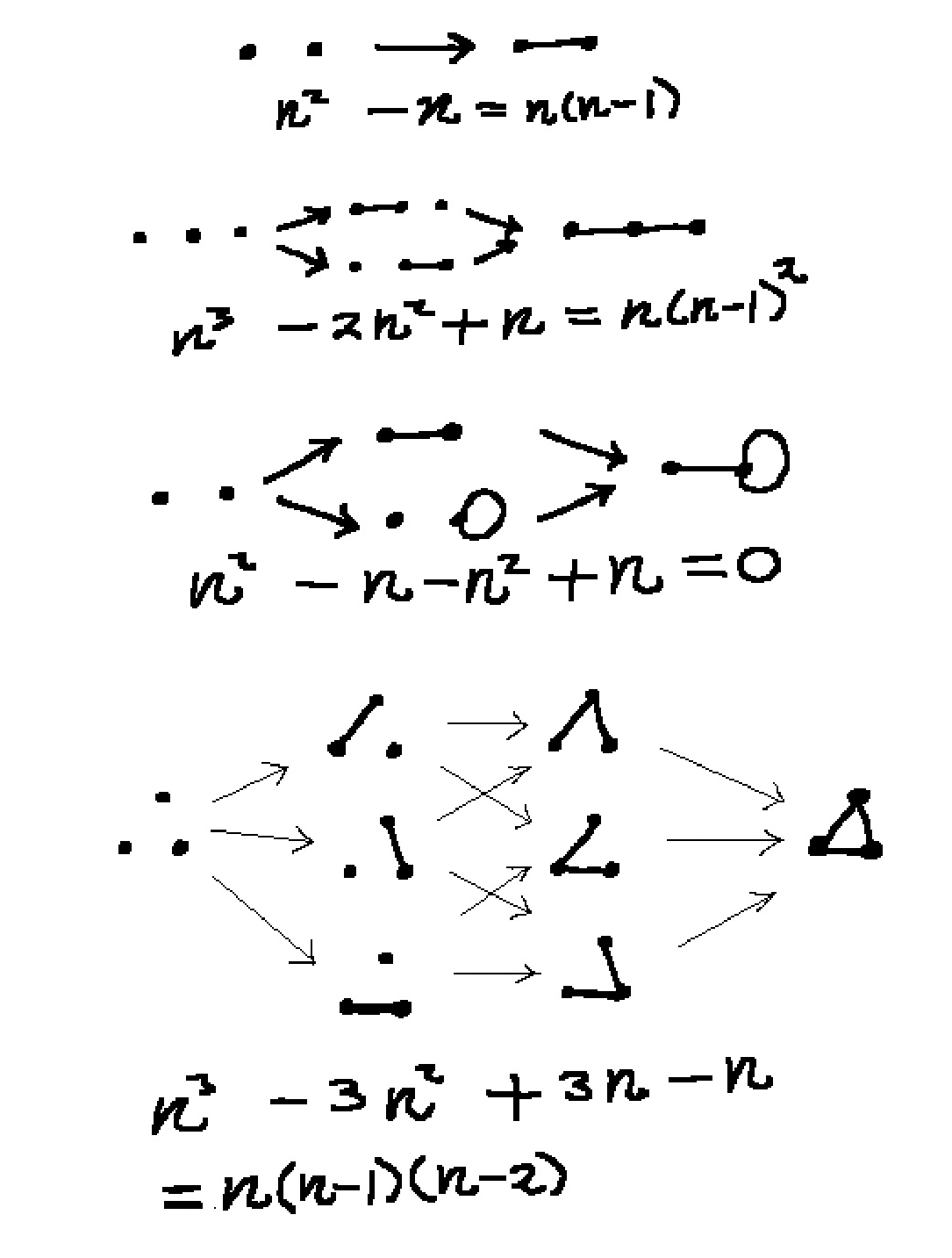}
\caption{Graph Categories}
\label{graphcomp}
\end{figure}

 In Figure~\ref{graphcomp} we illustrate how the subgraphs in this formula may be arranged from left to right, starting from the disjoint collection of the nodes of $G$ on the left, and ending in the full graph $G$ on the right. The vertical collections in the figure are subgraphs with the same number of edges. The arrows in the figure go from a subgraph with $k$ edges to a subgraph with $k+1$ edges. These collections of subgraphs are labeled with the contributions $\pm n^{c(S)}$ where $c(S)$ is the number of components of $S.$ The alternating sums yield the coloring numbers $C_{G}(n).$ This figure contains the seeds for the constructions in this paper. The arrow arrangement of subgraphs of $G$ can be regarded as a category with objects the subgraphs and arrows the generators of the non-identity morphisms of the category. This category then can be measured by choices of homology, using different aspects of algebra. The alternating sum becomes an Euler characteristic of the homology. The simplest version of this construction is given in the last section of this paper. The more complex versions involve chromatic polynomials.\\

The {\it color count} or {\it color evaluation} becomes the {\it chromatic polynomial} by replacing the integer $n$ by a polynomial variable $\lambda$ and then one has the chromatic polynomial as a sum over evaluations of subgraphs.
$$C_{G}(\lambda) = \sum_{S} (-1)^{e(S) }\lambda^{c(S)}$$ In order to categorify the chromatic polynomial we want the above state sum formula over subgraphs to be seen as a graded Euler characteristic. A strategy for accomplishing this is explained in Section 2. We change variables to  $\lambda =(1+q)$ and interpret $(1+q)^{c(S)}$ as a sum of contributions from {\it enhanced states} $s$ that consist in labelings of the components of $S$
with one of the two symbols $1$ or $x.$ If the enhanced state $s$ has $k$ labels $x$ then it contributes $q^k$ in the binomial expansion of  $(1+q)^{c(S)}.$  This allows us to rewrite the state sum as
$$C_{G}(q) = \sum_{S} (-1)^{e(S)} (1+q)^{c(S)} = \sum_{s} (-1)^{e(s)} q^{j(s)}$$ where $j(s)$ denotes the number of components of $s$ receiving the label $x.$ We take these labeled states (enhanced states) as generators of the chain complex where the $k$-th level $C^{k}(G)$ is the set of enhanced states with $k$ edges. The labels $x$ and $1$ are regarded as the generators of the ring $Z[x]/(x^2)$ and we define $d: C^{k}(G) \longrightarrow C^{k+1}(G)$ by (details in the next section) summing over the results of single edge additions to the subgraph. Each edge, when added, may go between components of a state, or it may connect a component with itself. This results in either two components merging to one component, or an addition of an edge to a given component. In the first case we assign the product of the labels on the two distinct components that are connected by the edge to be the label of the new component. In the case of a single component that remains a single component we retain its label. From this definition we see that for each addition of one edge to an enhanced state $s$ to form a new enhanced state $s'$ it is the case that $j(s) = j(s').$ That is, the number of $x$ labels does not change. Thus the boundary map preserves $j(s)$ and we have $d:C^{k,j}(G) \longrightarrow C^{k+1,j}$ where the $j$ denotes the subcomplex of chains generated by enhanced states $s$ with $j = j(s).$ It then follows that the chromatic polynomial is a graded Euler characteristic of the homology of this complex, and we have
$$\sum_{j} q^{j} \chi(H^{\star,j}(G))=\sum_{j} q^{j} \chi(C^{\star,j}(G)) = \sum_{s} (-1)^{e(s)} q^{j(s)} = C_{G}(\lambda)$$ where $\lambda = 1 + q$. This construction categorifies the chromatic polynomial, producing the same categorification as the one given in \cite{Rong1}. We explain the method here in some detail because we use the method of enhanced states uniformly in this paper for the new constructions that we will produce.\\

The key point about using enhanced states for categorification is that it can sometimes be arranged that the powers of a polynomial variable correspond to counts in the labelings of states in a state sum. Thus in the chromatic categorification above, we
have $q^{j}$ corresponding to states with $j$ labels of the form $x.$ And the algebra of $x$ and $1$ in the ring $Z[x]/(x^2)$ combined with the rules for the differential guarantee that the $j$-count is not changed under the differential. This divides the chain complex into parts with specific $j$ values and allows that the coefficient of $q^{j}$ in the polynomial can be an Euler characteristic of the homology at grading $j.$\\

Section 3 gives a categorification of the dichromatic polynomial, which is defined by the equations below (see the section 3 for details).
$$Z_{G}(v, \lambda) = Z_{G - e}(v, \lambda) + v Z_{G/e}(v, \lambda)$$
and
$$Z_{G \sqcup \bullet}(v, \lambda) = \lambda Z_{G}(v, \lambda).$$
The variable $v$ in the dichromatic polynomial generalizes the minus one in the chromatic polynomial. We use an enhanced states formulation with $\lambda = 1 + q$ and $v = 1 + p,$ and find a natural triply graded homology theory for these variables so that 
$$Z_{G}(p,q) =\sum_{i,j} p^{i}q^{j} \chi (H^{\star, i,j}(G)).$$\\

Furthermore, we define a sequence of generalizations of the chromatic polynomial that we dub {\it impropriety polynomials} $C^{i}_{G}(\lambda)$ that count the number of {\it improper} colorings of a graph $G$ where an {\it impropriety}
is an edge whose endpoints receive the same color or a loop whose endpoint (must) receives the same color. The $i$-th impropriety polynomial counts the number of colorings of $G$ with exactly $i$ improprieties.
These new polynomials are related to the dichromatic polynomial as follows. If we let $v = \alpha -1$ so that 
$$Z_{G}(\alpha, \lambda) = Z_{G - e}(\alpha, \lambda) + (\alpha -1) Z_{G/e}(\alpha, \lambda)$$
then 
$$Z_{G}(\alpha, \lambda) = \sum_{i=0}^{\infty} \alpha^{i} C^{i}_{G}(\lambda).$$ 
Thus with this change of variables, the dichromatic polynomial can be regarded as a series in the variable $\alpha$ whose coefficients are the impropriety polynomials.\\

This result has the following striking corollary, that has application to the Potts model in statistical mechanics.\\

\noindent{\bf Corollary.} In the variable $\alpha$ as above, the dichromatic polynomial for an integer number of colors can be expressed as follows $$Z_{G}(\alpha, n) = \sum_{\sigma} \prod_{\langle i,j \rangle} \alpha^{\delta(\sigma(i), \sigma(j))}$$
 where $\sigma: Nodes(G) \longrightarrow \{1,2,3,...,n\}$ runs over all color assignments to the nodes of 
$G,$ $\langle i,j \rangle$ denotes an edge of $G$ with end nodes $i$ and $j$ and $\delta(a,b)$ is the Kronecker delta, equal to $0$ when $a \ne b$ and equal to $1$ when $a = b.$\\
Thus $$\sum_{i=0}^{\infty} \alpha^{i} C^{i}_{G}(n) = \sum_{\sigma} \prod_{\langle i,j \rangle} \alpha^{\delta(\sigma(i), \sigma(j))}.$$\\

We find that the impropriety polynomials are themselves graded Euler characteristics in the homology theory of the dichromatic polynomial.
$$C^{i}_{G}(\lambda) = (-1)^{i} \sum_{j} (\lambda-1)^{j} \chi (H^{\star, i, j}(G)).$$ \\

Section 4 discusses the Potts model in statistical mechanics in the light of the results in section 3. We refer the reader to this section for the definition of the Potts model and its properties.
A main result of this section has the form 
$$Z_{G}(T) = Z_{G}(p,q) =\sum_{i,j} (-1)^{i} (e^{\frac{-1}{kT}})^{i}(n-1)^{j} \chi (H^{\star, i,j}(G)).$$
It is worth remarking that the terms $(e^{\frac{-1}{kT}})^{i}$ correspond to specific energy levels in the partition function for the Potts model. Thus the coefficients of these energy levels are expressed
in terms of graded Euler characteristics in the dichromatic homology corresponding to level $i$ impropriety. In the Potts model, high impropriety corresponds to high energy and proper coloring is the ground state at zero energy.\\

Section 5 uses the enhanced state methodology to categorify the Penrose-Kauffman polynomial 
\cite{P, VCP,LK1,LK2,LK3,LK4,BKM,BMC} which is defined for special $n$-colorings of trivalent graphs $G$ with specified perfect matchings. 
$$PK[ \YMGlyph ] = PK[ \CrossStarGlyph ] = PK[ \VDiag \ + PK[ \CDiag ] - 2PK[ \CDotDiag ],$$
$$PK[ G O ] = \lambda PK[G],$$
$$PK[ O ] = \lambda.$$\\
See the section 5 for the results of categorification and see the references above for background on these polynomials.\\

Section 6 returns to the chromatic evaluation at a choice of $n$-colors. We have the numbers $C_{G}(n)$ and we show that there are many ways to choose an algebra ${ A}$ and define differentials for a homology theory so that these evaluations are Euler characteristics of a homology theory. Thus $$C_{G}(n) = \chi(H^{\star}(G, { A})).$$ We give a specific family of algebras ${ A}(n)$ so that proper colorings correspond one-one to 
homology elements of grading $0$ and so that all the higher homology groups vanish. There are other choices of algebra that are related to specific coloring problems. For example the integral group ring of the Klein 4-group is an approriate algebra for the 4-coloring problem. It is planned to return to all the parts of this paper with more analysis and more calculations. We are particularly interested to see the relationships between coloring algebras in the sense of this last section and problems for the proper colorings of graphs.\\

\noindent {\bf Remark.} There have been a number of papers on homology related to chromatic polynomials. We give the following list in our bibliography
\cite{Luse,Rong1,Rong2,Stosic1,Stosic2,Jozef,LKPotts,Radmila,Bannerjee,BMC,BMV}. The present paper is intended to be self-contained, and we have made reference to this literature
when it is appropriate. Many ideas in these papers and in the present paper come from the structure of Khovanov homology for knots and links. We do not refer directly to that literature but point out that the papers by 
Bar-Natan, Turner and the author \cite{BarNatan,Turner,LKNew,LKho1,LKho2} are a useful linkage with categories, combinatorics and topology.\\

\section{Categorifying Chromatic Polynomials}
Recall the chromatic polynomial $C_{G}(\lambda)$ defined as a polynomial in the ring  $Z[\lambda]$ for any finite graph $G.$
The polynomial is defined so that for any natural number $n$, $C_{G}(n)$ is equal to the number of proper vertex colorings of the graph $G$ using
$n$ colors. An assignment of colors to the vertices of a graph is said to be {\it proper} if distinct colors are assigned to any two vertices that are connected by an edge in the graph.
Hassler Whitney \cite{Whitney} observed that, given an edge $e$ in $G,$ then $$C_{G}(n) = C_{G - e}(n) - C_{G/e}(n)$$ where $G-e$ denotes the graph obtained by deleting the edge $e$ and $G/e$ denotes the graph obtained by 
removing the edge $e$ and identifying the end nodes of $e$ with each other. This identity is a tautology when the graph $G$ is {\it simple}, meaning that two nodes either have one edge or no edge between them. It remains true for graphs with multiple edges. Along with Whitney's logical identity (above) we also have that $$C_{G \sqcup \bullet}(n) = n C_{G}(n)$$ where $\bullet$ denotes an isolated node. These two properties completely deternine 
$C_{G}(n)$ and one can use them to see that there is the formula $$C_{G}(n) = \sum_{S} (-1)^{e(S)}n^{c(S)}$$ where S runs over all subgraphs of $G$ (from a collection of disjoint nodes to the full graph $G$) and 
$e(S)$ is the number of edges in $S$ and $c(S)$ is the number of connected components in $S.$\\

\noindent {\bf Remark.} One way to prove this summation formula is to use the Whitney identity as a recursion and mark each edge that is to be contracted. Then the recursion yields a sum over states of the graph $G$ where every edge has either been deleted or marked. These states correspond to all the subgraphs of $G$. Each component subgraph is to be contracted to one node and so contributes $n$ possible colors. The sign of the contribution of a given subgraph
corresponds to the parity of the number of marked edges, and hence to the parity of the number of edges of the subgraph.\\

\noindent {\bf Remark.} Another way to prove the chromatic summation formula is to regard it as an application of the inclusion-exclusion principle. Let a color assignment to a graph be an assignment of colors to each of the components of the graph. One counts all color assignments to the nodes of the graph, then subtracts all color assignments to subgraphs with one edge, then adds the count of all color assignments to subraphs with two edges and so on. The resulting sum constitutes all colorings where nodes connected by edges are assigned distinct colors.\\

The chromatic polynomial is obtained by replacing the $n$ in the chromatic evaluation $C_{G}(n)$ with a polynomial variable $\lambda$. Thus we have
$$C_{G}(\lambda) = \sum_{S} (-1)^{e(S)} \lambda^{c(S)}$$
and the consequences
$$C_{G}(\lambda) = C_{G - e}(\lambda) - C_{G/e}(\lambda)$$
and
$$C_{G \sqcup \bullet}(\lambda) = \lambda C_{G}(\lambda).$$\\

\noindent{\bf Categorification of the Chromatic Polynomial.}
We begin by pointing out how one can {\em categorify} the chromatic polynomial so that it is seen as a graded Euler characteristic of a homology theory associated with the graph.
Call the subgraphs $S$ in the summation formula for the chromatic polynomial the {\it states} of $G.$  Take a variables $x$ and $1$ from the ring $M=Z[x]/(x^2)$ as generators for $M$ as a module over $Z.$
We will use $x$ and $1$ as labels and then bring in the algebra structure in $M.$ Call a state with labels of either $x$ or $1$ on its {\it components} an enhanced state $s$ (we will use lower case $s$ for enhanced states).
These algebraic labels will be used to define generators and differentials for a chain complex. Let $\lambda = (1 + q)$ where $q$ is another polynomial variable. Note that there is a $1-1$ correspondence between labelling $c(S)$ components with $1$ or $x$ and products using $1$ and $q$ for each component. For example if $c(S) = 2$ then
we would have $$11,1x,x1,xx$$ as labelings of two (ordered) components, and $$(1+q)^2 = 1\times1 + 1 \times q + q \times 1 + q \times q.$$ Given an enhanced state, we let $j(s)$ denote the number of 
components of that state that have a label equal to $x.$ Thus we can regard $$\lambda^{c(S)} = (1+q)^{c(S)}$$ in the state summation
as corresponding to a sum over all the enhancements of the state $S$, counting by $q^{j(s)}$  where $j(s)$ is the number of $x$-labels on the enhanced state. Thus we can rewrite the chromatic state sum in terms of enhanced states as shown below.
$$C_{G}(q) = \sum_{S} (-1)^{e(S)} (1+q)^{c(S)} = \sum_{s} (-1)^{e(s)} q^{j(s)}$$
Here the number of edges $e(s)$ of an enhanced state $s$ is the same as the number of edges of the state $S.$\\

The chain complex ${ C}(G),$  associated with $G,$ is generated by the enhanced states of $G$. A generator $s$ of ${ C}(G)$ is a subgraph of $G$ where each component of the subgraph is labeled with 
either $x$ or $1.$ Let ${ C}^{k}(G)$ denote that part of the complex that is generated by states with $k$ edges. We will define a differential $d$ for the chain complex so that 
$$d: { C}^{k}(G) \longrightarrow { C}^{k+1}(G).$$ First we define partial differentials $\partial_{e}(s)$ that transform an enhanced state with $k$ edges to a new enhanced state with $k+1$ edges.
Then $d$ is a sum of partial differentials as described below. To define the partial differential, choose an edge $e$ in the graph $G$ that is not in the state $s$. (If there is no such edge, then the differential will be zero.) 
The ends of the edge $e$ are in either one or two labeled components of $s$. \\

Let s[e] denote the new state obtained by adding the edge $e$ to the state $s$, making a new subgraph with one more edge.
Let $s[e, a]$ denote this new state with label $a$ attached to the component of $s[e]$ that contains $e$. We will set the value of this label with the partial differential.\\

If there are two components with labels $a$ and $b$ corresponding to the ends of $e$, define $$\partial_{e}(s) = s[e,ab].$$ If $ab=0$ then we take $s[e,0] = 0$ so that the partial differential vanishes in this case.
Thus, if $a = b = 1$ then $\partial_{e}(s) = s[e,1].$
If $a = 1$ and $b = x$ or  $a = x$ and $b = 1$ then $\partial_{e}(s) = s[e,x].$ If $a = x$ and $b = x$  then $\partial_{e}(s) = 0.$ \\
 
If the the ends of $e$ are in a single component with label $a$, then $\partial_{e}(s) = s[e,a].$ 
It is important to note that the number of $x$ labels remains the same under the action of the partial differentials.\\

To define the total differential modulo two, we sum over all edges in the $G$ complement of the state $s.$ Thus, modulo two we have $$d(s) = \sum_{e \in G-s} \partial_{e}(s).$$
Here by $e \in G-s$ we mean that $e$ is an edge in $G$ that is not in $s.$ It is not hard to see that $d(d(s)) = 0$ modulo two since the partial differentials commute with each other.
The full differential over the integers is obtained by {\it ordering the edges of $G$} and letting $$d(s) = \sum_{e \in G-s} (-1)^{n(e)} \partial_{e}(s),$$ where $n(e)$ denotes the number of edges in $s$ that precede $e$
in the edge-order for $G.$ One can show that $d^2 = 0$ and that the homology groups obtained from the differential are independent of this choice of order.\\

We can now define $H^{k}(G)$ as the kernel of $d: { C}^{k}(G) \longrightarrow { C}^{k+1}(G)$ modulo the image of $d: { C}^{k-1}(G) \longrightarrow { C}^{k}(G).$ 
In order to see how these homology groups of the graph $G$ are related to its chromatic polynomial, note that $j(s)$, the number of $x$ labels in the enhanced state $s$, is not changed under the application of the differential.
Thus we can refine this homology with a second grading by defining ${ C}^{k,j}(G)$ to be the submodule of chains that have a fixed number of edges $k$ and a fixed number of $x$ labels $j.$ Then
$$d: { C}^{k,j}(G) \longrightarrow { C}^{k+1,j}(G)$$ and we have graded homology groups $H^{k,j}(G).$ \\

Now return to the formula for the chromatic polynomial.
$$C_{G}(q) = \sum_{s} (-1)^{e(s)} q^{j(s)} = \sum_{j} q^{j} \sum_{s: j(s)=j} (-1)^{e(s)}$$
$$ = \sum_{j} q^{j} \sum_{e} (-1)^{e} rank{ C}^{e,j}$$
$$ = \sum_{j} q^{j} \sum_{e} (-1)^{e} rank H^{e,j}$$
$$C_{G}(q) =\sum_{j} q^{j} \chi (H^{\star, j}(G)).$$
In other words, we see that the original formula, counting with signs the numbers of enhanced states with a given edge count and label count can be re-arranged as a sum over powers of 
the variable $q$ corresponding to label counts multiplied by Euler characteristics of the chain complexes in the $j$-grading. This becomes Euler characteristics of the homology by the usual relationship between homology and 
chain complex. This means that he coefficients of the chromatic polynomial $C_{G}(q)$ in the variable $q$ are Euler characteristics in this homology theory. It still remains to be seen what this result says about the chromatic polynomial. It is not hard to see that the homology theory itself contains more information about the graph than just the chromatic polyomial.\\

The key point about this categorification is that it arises from looking at a refined structure of subgraphs of a graph. By labeling the components of each subgraph we obtain a way to create a directed dynamic letting given subgraphs associate to other subgraphs by joining along edges in the complement of the subgraph. Just so in statistical mechanics one has the possibility of regions of constant spin merging with one another. There is a close analogy with such dynamics and the dynamics of the algebraic compositions in the graph homology. The combinatorial fact that the differential preserves labeling counts on regions of constant spin suggests further relationships with the statistical physics.\\

It is our purpose in this paper to show how this kind of categorification of the chromatic polynomial extends to the dichromatic polynomial and how this extension is related to the Potts model in statistical mechanics.
We will also show how these methods can be used to categorify the Penrose-Kauffman polynomial \cite{P,BKM,LK4} and a dichromatic extension of that polynomial. Finally, we will examine variations on these structures by categorifying coloring evaluations rather than polynomials.\\

\section{The Dichromatic Polynomial}
The dichromatic polynomial $Z_{G}(v,\lambda)$ of a graph $G$ is defined by adding a new variable $v$ to the chromatic polynomial. The defining relations are given by the equations below, where the summation is over all subgraphs of $G$, $e(S)$ denotes the number of edges in the subgraph $S$ and $c(S)$ denotes the number of components in the subgraph $S.$
$$Z_{G}(v, \lambda) = \sum_{S} v^{e(S)} \lambda^{c(S)}$$
and the consequences
$$Z_{G}(v, \lambda) = Z_{G - e}(v, \lambda) + v Z_{G/e}(v, \lambda)$$
and
$$Z_{G \sqcup \bullet}(v, \lambda) = \lambda Z_{G}(v, \lambda).$$\\

Let $\lambda = 1 + q$ and $v = 1 + p.$ Choose labels $1$ and $y$ for the edges of subgraphs. Define a {\it refined enhanced state s} of a subgraph $S$ of $G$ to be {\it a labeling of the components of 
$S$ by either $1$ or $x,$ and an independent labeling of the edges of $S$ with either $1$ or $y$.} The labels $1$ and $y$ for edges are indicative and not elements of an algebra. They will be used to extend the differential defined in the first section of the paper. \\

We now have
$$Z_{G}(v,\lambda) = Z_{G}(p,q) = \sum_{S} (-1)^{e(S)}(1+p)^{e(S)} (1+q)^{c(S)} = \sum_{s} (-1)^{e(s)} p^{i(s)}q^{j(s)}$$
where the last sum is over all refined enhanced states $s$, $e(s)$ is the number of edges in $s,$ $i(s)$ is the number of $y$ labels on the edges of $s$ and $j(s)$ is the number of $x$ labels on the components of $s.$ \\

\noindent{\bf Dichromatic Chain Complex.}
We define a new chain complex ${ C}(G)$ generated by the refined enhanced states of $G.$ This chain complex is bigraded where ${ C}^{e,i,j}(G)$ is generated by enhanced states $s$ with $e$ edges, and
$i = i(s)$, the number of edges in $s$ labeled with $y$, and $j = j(s)$, the number of components of $s$ with label $x.$ The partial differentials are extended.  The partial differentials add new edges and assign $1$ to the new edge, and assign $1$ or $x$ to the component containing the new edge, using the algebraic rules described in the last section via $x^2=0$ and $1$ as the identity in the algebra.\\

This extension of the differential ensures that both $i$ and $j$ are preserved under the partial differentials. Thus we now have a chain complex ${ C}^{\star,i,j}(G)$ with $d: { C}^{k,i,j}(G) \longrightarrow { C}^{k+1,i,j}(G)$  and bigraded homology groups $H^{k,i,j}(G).$ The dichromatic polynomial is a graded Euler characteristic in this bigraded theory. Thus we have\\

\noindent{\bf Theorem.} $$Z_{G}(p,q) =\sum_{i,j} p^{i}q^{j} \chi (H^{\star, i,j}(G)).$$\\

\begin{figure}
\includegraphics[width=9cm]{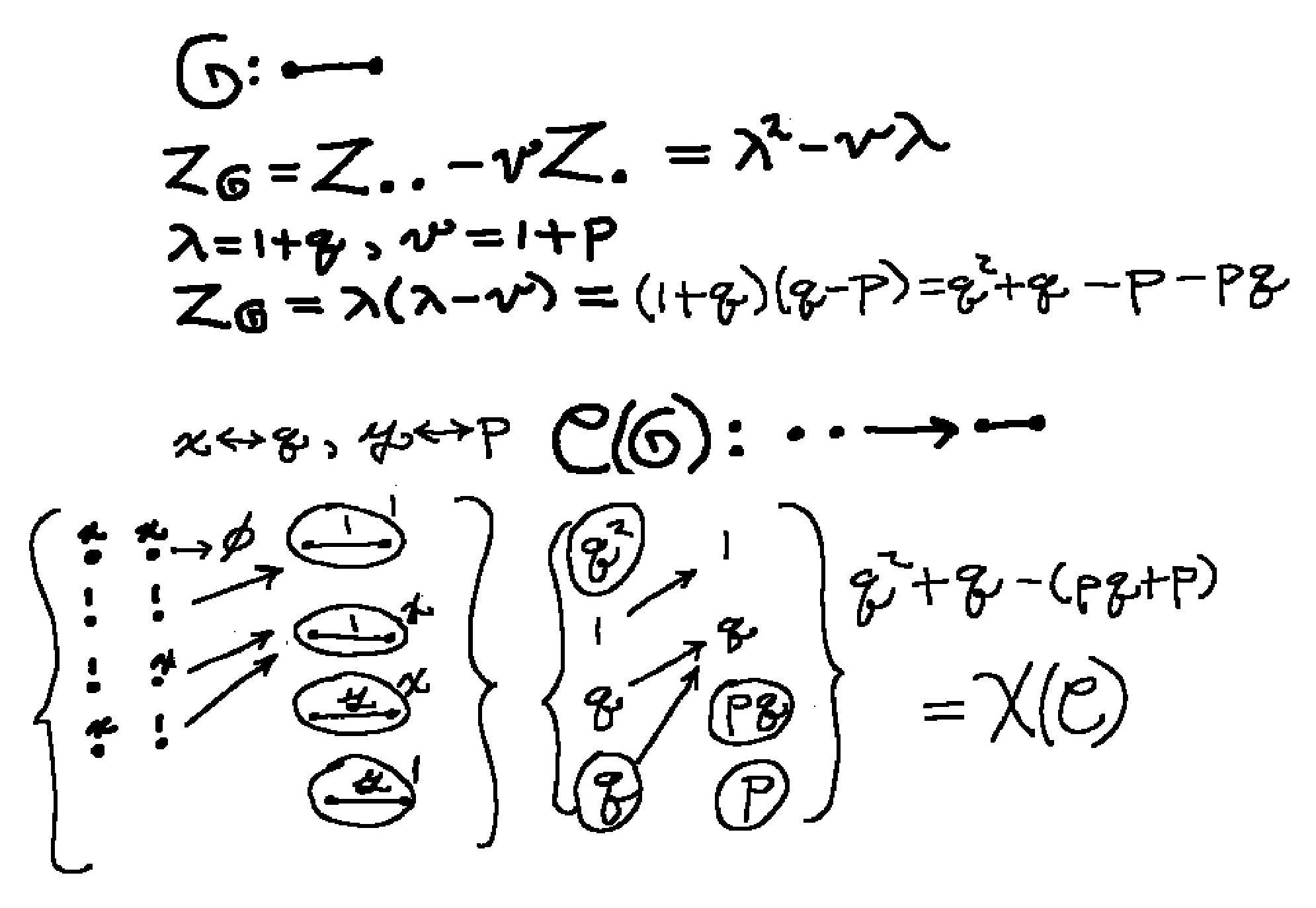}
\caption{\bf Dichromatic Homology for a simple line graph.}
\label{dichrome}
\end{figure}

In Figure~\ref{dichrome} we illustrate the dichromatic polynomial and its chain complex in the case where the Graph $G$ has two distinct nodes and a single edge. In this case $Z_{G} = \lambda^2 - \lambda v = \lambda (\lambda - v)$ from the definition of the dichromatic polynomial. Substituting $\lambda = 1+ q$ and $v = 1 + p$ we find that $Z_{G} = (1+q)(q-p) = q^2 + q - p - pq.$ The figure illustrates how these terms correspond to the graded Euler characteristic 
of the refined enhanced states that generate the chain complex. One can also see how the homology of the complex works out so that 
$$H^{0,2,0}=Z, H^{0,0,0} = 0, H^{0,1,0}=Z,H^{1,0}=0, H^{1,0,0}=0, H^{1,1,0}=0,H^{1,0,1}=Z,H^{1,1,1}=Z.$$
As in the figure we can write the graded Euler characteristic of the complex ${ C}(G)$ as the dichromatic polynomial.
$$\chi({ C}(G)) = q^2 + q - p - pq.$$\\

\begin{figure}
\includegraphics[width=9cm]{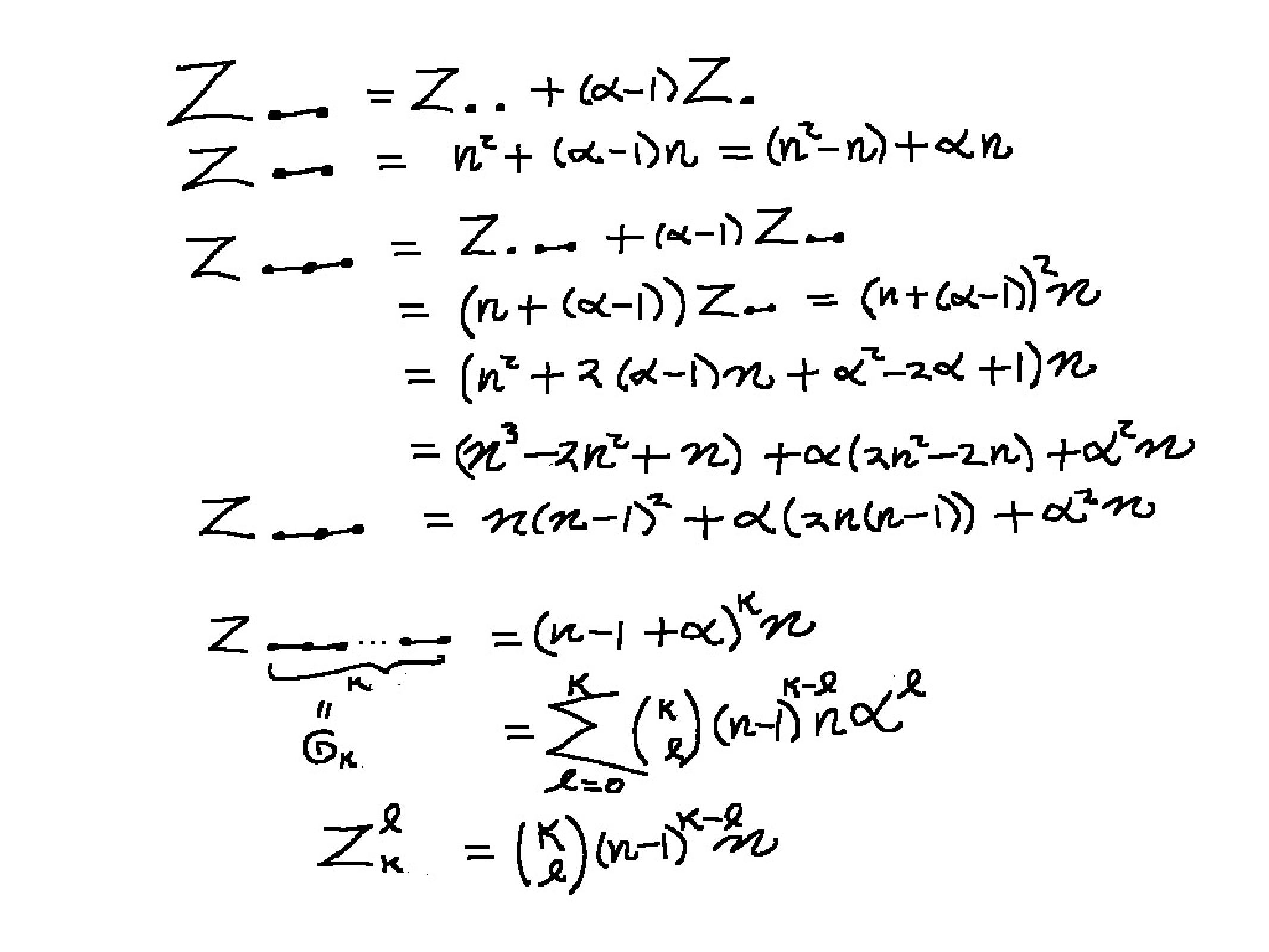}
\caption{\bf Dichromatic Examples}
\label{dichromeexamples}
\end{figure}

\noindent{\bf Remark.} Consider the dichromatic polynomial with $v = 1 - \alpha$, that is, with $p = -\alpha.$ Then the basic identities for the polynomial are as follows, taking $\lambda = n$ so that we can discuss
what is counted by this polynomial:
$$Z_{G}(\alpha, n) = Z_{G - e}(\alpha, n) + (\alpha -1) Z_{G/e}(\alpha, n)$$
and
$$Z_{G \sqcup \bullet}(\alpha, n) = n Z_{G}(\alpha, n).$$\\

\noindent{\bf Theorem.} $Z_{G}(\alpha, n) = \sum_{i=0}^{\infty} \alpha^{i} C^{i}_{G}(n)$ where $C^{i}_{G}(n)$ is the number of colorings of the edges of $G$ with $n$ colors so that {\it the endpoints of exactly $i$ edges receive the
same colors.} $C^{i}_{G}(n)$ counts the number of colorings of $G$ with {\it impropriety $i$}. Thus $C^{0}_{G}(n) = C_{G}(n)$ is the usual chromatic evaluation of $G.$ By replacing $n$ by an algebraic variable $\lambda$ each polynomial coefficient $C^{i}(\alpha, \lambda)$ can be taken as a generalization
of the chromatic polynomial $C^{0}_{G}(\lambda) = C_{G}(\lambda)$ for impropriety $i.$\\

\noindent{\bf Proof.} Start with $Z_{G}(\alpha, n) = \sum_{i=0}^{\infty} \alpha^{i} C^{i}_{G}(n)$ as a definition of $Z_{G}(\alpha, n)$ with $C^{i}_{G}(n)$ defined as the count of colorings of $G$ that are proper except for exactly 
$i$ edges. Then note that the equation $Z_{G}(\alpha, n) = Z_{G - e}(\alpha, n) + (\alpha -1) Z_{G/e}(\alpha, n)$ is a logical identity about $Z_{G}(\alpha, n).$ To wit, any coloring of the ends of an edge $e$ is either proper. in which case there will be no contribution from $G/e$, or the two ends of $e$ receive the same color. In that case the $-1$ in  $(\alpha -1) Z_{G/e}(\alpha, n)$ cancels the color count from $Z_{G - e}(\alpha, n)$ but the 
$\alpha$ term marks and multiplies this coloring of $G/e.$ The upshot is that the logical identity ensures that colorings with $i$ improper edges will be counted with the coefficient $\alpha^{i}.$ This is what we needed to prove.
$\hfill\Box$\\ 

\noindent{\bf Example.} Let $G$ be the graph with two nodes and one edge. Then $Z_{G}(\alpha, n) = Z_{\bullet \,\, \bullet} + (\alpha-1) Z_{\bullet}$ $ = n^2 + (\alpha-1)n = n(n-1) + \alpha n.$
Indeed, $C_{G}(n) = C^{0}_{G}(n) = n(n-1)$ and $C^{1}_{G}(n) = n.$ In Figure~\ref{dichromeexamples} we illustrate the generalization of this formula for line graphs with $k$ edges. Letting $G_{k}$ denote such a graph, the figure illustrates the computation that shows $$Z_{G_{k}} = \sum_{l=0}^{k} {k\choose l} (n-1)^{k-l}n \alpha^{i}.$$ Thus $$C^{l}_{k}(n) = {k\choose l} (n-1)^{k-l}n$$ and it is easy to independently verify that this is the case since
choosing $l$ improper edges from the $k$ edges yields the proper coloring problem for a straight line graph with $(k-l)$ edges.\\ 

\noindent{\bf Corollary.} $Z_{G}(\alpha, n) = \sum_{\sigma} \prod_{\langle i,j \rangle} \alpha^{\delta(\sigma(i), \sigma(j))}$ where $\sigma: Nodes(G) \longrightarrow \{1,2,3,...,n\}$ runs over all color assignments to the nodes of 
$G,$ $\langle i,j \rangle$ denotes an edge of $G$ with end nodes $i$ and $j$ and $\delta(a,b)$ is the Kronecker delta, equal to $0$ when $a \ne b$ and equal to $1$ when $a = b.$\\

\noindent{\bf Proof.} It follows immediately that $\sum_{\sigma} \prod_{\langle i,j \rangle} \alpha^{\delta(\sigma(i), \sigma(j))} = \sum_{i=0}^{\infty} \alpha^{i} C^{i}_{G}(n),$ and so the result follows from the Theorem above.
$\hfill\Box$\\ 

\noindent{\bf Remark.} Note that it follows from the lemma that $C_{G}(n) = Z_{G}(0, n) = \sum_{\sigma} \prod_{\langle i,j \rangle} 0^{\delta(\sigma(i), \sigma(j))}.$  In this formula we count colorings via the properties of zero: $0^{0} =1, 0^{1} = 0.$\\

Finally, note the conseqence for evaluation of  the improprieties $C^{i}_{G}(\lambda)$ in terms of the dichromatic homology. We have that $Z_{G}(\alpha, \lambda)$ is obtained from $Z_{G}(p, q)$  with $p= -\alpha$ and $q = \lambda-1.$ We have shown that $Z_{G}(p,q) =\sum_{i,j} p^{i}q^{j} \chi (H^{\star, i,j}(G))$ and so we have
$$Z_{G}(\alpha, \lambda) =\sum_{i,j} (-1)^{i} \alpha^{i}(\lambda-1)^{j} \chi (H^{\star, i,j}(G)) = \sum_{i} \alpha^{i} [(-1)^{i} \sum_{j} (\lambda-1)^{j} \chi (H^{\star, i,j}(G))].$$
Thus \\

\noindent{\bf Theorem.} The impropriety polynomials are expressed in terms of graded Euler characteristics of the homology.
$$C^{i}_{G}(\lambda) = (-1)^{i} \sum_{j} (\lambda-1)^{j} \chi (H^{\star, i, j}(G)).$$  \\

In the next section we apply these results to the Potts model in statistical mechanics.\\

\section{The Potts Model}
Given a graph $G,$ the {\it partition function} for the Potts model with $n$ spin states at the nodes of the graph is given by the formula
$$Z_{G}(T) =\sum_{\sigma} e^{\frac{-E(\sigma)}{kT}}$$ 
Here $T$ is the temperature.  $\sigma$ denotes a {\it physical state} $\sigma : Nodes(G) \longrightarrow \{1,2,..., n \}$ as an arbitrary assignment 
of spins to the nodes of the graph. The {\it energy} $E(\sigma)$ is given by the formula $$E(\sigma) = \sum_{\langle i,j \rangle} \delta(\sigma(i), \sigma(j))$$
where $\langle i,j \rangle$ runs over all edges of $G$ with end nodes indicated by $i$ and $j,$ and $\delta(r,s)$ denotes the Kronecker delta that is equal to $1$ where $r=s$ and is equal to $0$ when
$r \ne s.$\\

The partition function for a statistical mechanics model contains a wealth of information, and as a total summation it can be interpreted to give the probability of the system being in a given state
at a given temperature $T$  via the formula
$$Prob(\sigma) =  e^{\frac{-E(\sigma)}{kT}}/Z_{G}(T).$$

It was observed by Neville Temperley \cite{Baxter,Temperley,LKStat} that the Potts model partition function is a specialization of the dichromatic polynomial of the graph $G.$\\

\noindent {\bf Theorem (Temperley).} Let $Z_{G}(T)$ be the Potts model partition function as described above. Then, using the definition of the dichromatic polynomial given in the last section, we have
that $Z_{G}(T) = Z_{G}(\lambda, v),$ the dichromatic polynomial for $G$ evaluated for $q$ spin states and $v = 1 - e^{\frac{-E(\sigma)}{kT}}.$ and $\lambda = n.$\\

\noindent {\bf Proof.} Note that we have 
$$Z_{G}(T)=\sum_{\sigma} e^{\frac{-E(\sigma)}{kT}}= \sum_{\sigma} e^{\frac{-\sum_{\langle i,j \rangle} \delta(\sigma(i), \sigma(j))}{kT}}$$
$$= \sum_{\sigma} \prod_{\langle i,j \rangle} (e^{\frac{-1}{kT}})^{\delta(\sigma(i), \sigma(j))}.$$ From this one sees, that if one fixes attention on a single edge in the graph, that with the edge absent there is an effective contribution of a factor of $1$ in the product in the last formula, and if the edge is present there is a factor of $e^{\frac{-1}{kT}}$ in the formula. Since $1 - (1 - e^{\frac{-1}{kT}}) = e^{\frac{-1}{kT}}$, we see that the formula
$$Z_{G} = Z_{G - e} - (1 - e^{\frac{-1}{kT}})  Z_{G/e}$$ follows from the definition of the Potts Partition function. The rest of the proof follows similarly. $\hfill\Box$\\

\noindent{\bf Remark.} Note that, letting $p$ replace $-e^{\frac{-1}{kT}},$ we can define 
$$Z_{G} = \sum_{\sigma} \prod_{\langle i,j \rangle} p^{\delta(\sigma(i), \sigma(j))}.$$
and deduce that 
$$Z_{G} = Z_{G - e} - (1 + p)  Z_{G/e},$$ matching our formalism for the dichromatic polynomial with $v = (1 + p).$
Note also that when $p =0$ we retrieve the chromatic polynomial at specific values of $\lambda = n$ with 
$$C_{G}(n) = \sum_{\sigma} \prod_{\langle i,j \rangle} 0^{\delta(\sigma(i), \sigma(j))}.$$
This is, in its way, a remarkable formula for the chromatic evaluation, based on the basic mathematical equations $ 0^{1} = 0$ and $0^{0}=1$ so that the formula only counts proper colorings of the graph $G.$ In the Potts model
the value $p=0$ corresponds to the zero-temperature limit. (Recall that $T$ in $e^{\frac{-1}{kT}}$ is the temperature variable in the model.) $$Lim_{T \longrightarrow 0} e^{\frac{-1}{kT}} = 0.$$\\

Using the dichromatic homology of the last section, we can now express the Potts partition function as a graded Euler characteristic.
We let $\lambda = n = 1 + q$ and $v = 1 + p.$ Thus $p = - e^{\frac{-E(\sigma)}{kT}}.$ We have the formula
$$Z_{G}(p,q) =\sum_{i,j} p^{i}q^{j} \chi (H^{\star, i,j}(G)).$$
and this translates to
$$Z_{G}(p,q) =\sum_{i,j} (-1)^{i} (e^{\frac{-1}{kT}})^{i}(n-1)^{j} \chi (H^{\star, i,j}(G)).$$\\
Thus we have proved the \\

\noindent {\bf Theorem.} Let $Z_{G}(T)$ denote the Potts model partition function for a graph $G$ with $n$ spin states, temperature variable T  and Boltzmann constant k.
Then $Z_{G}(T)$ is given as a graded Euler characteristic for the dichromatic bigraded homology of the graph, with specific formula
$$Z_{G}(T) = Z_{G}(p,q) =\sum_{i,j} (-1)^{i} (e^{\frac{-1}{kT}})^{i}(n-1)^{j} \chi (H^{\star, i,j}(G)).$$
Note that this formula gives the partition function in terms of powers of the number of spins less one ($n-1$) and powers of $e^{\frac{-1}{kT}}$ with coefficents the corresponding Euler characterisics 
of the dichromatic graded homologies.\\

\noindent {\bf Remark.} In further work we will analyze the physics of the Potts model in terms of properties of the homology. In particular, the coefficents of $(e^{\frac{-1}{kT}})^{i}$ (corresponding to specific energy levels) are the impropriety polynomials
defined in the previous section and are described homologically. We will explore the relationships of the dichromatic homology and the properties of the statistical mechanics model - phase transition and criticality.\\

\begin{figure}
    
     \includegraphics[width=9cm]{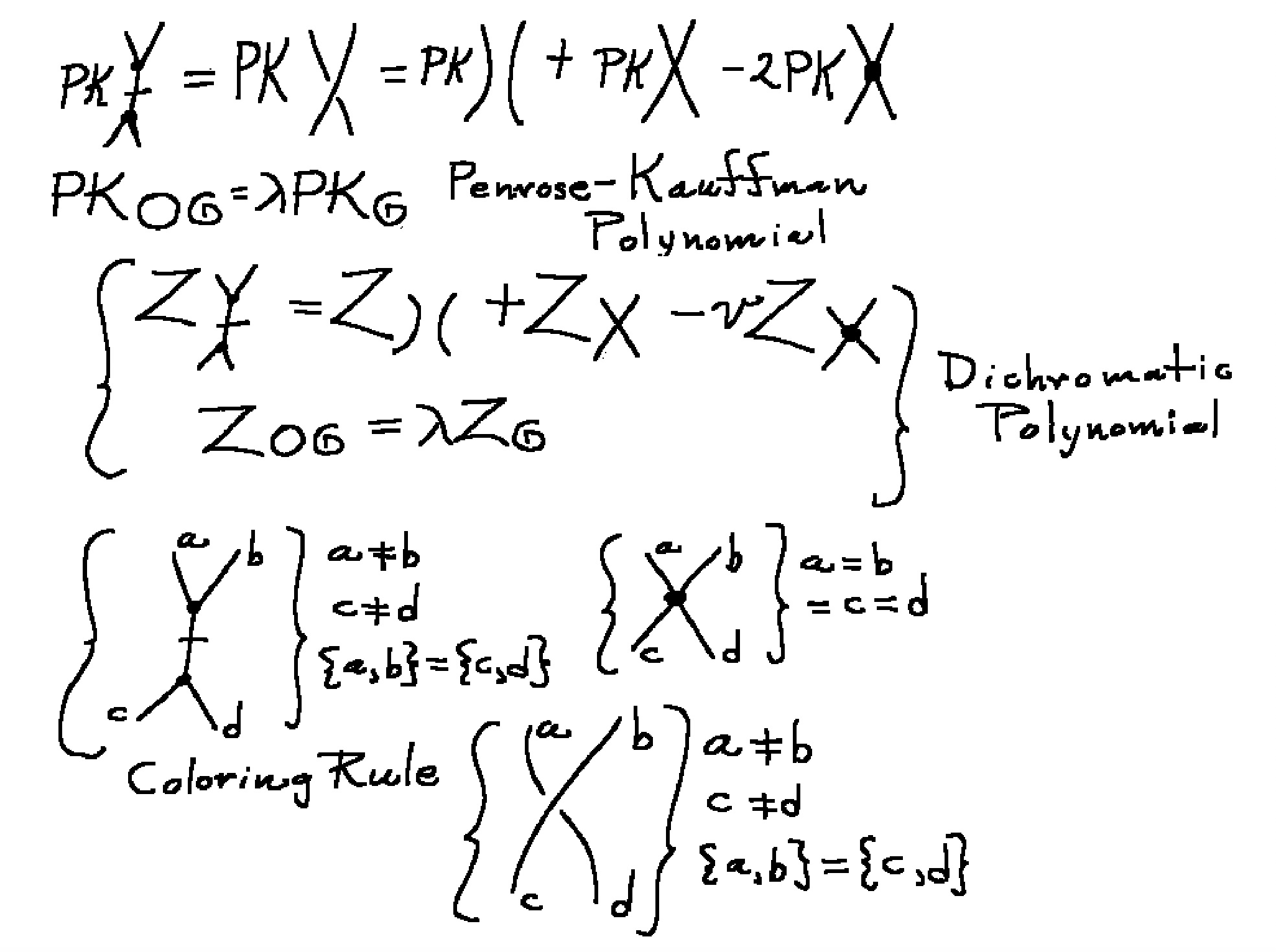}
    
     \caption{\bf Penrose-Kauffman Polynomials}
     \label{penrose}

\end{figure}

\begin{figure}
    
     \includegraphics[width=9cm]{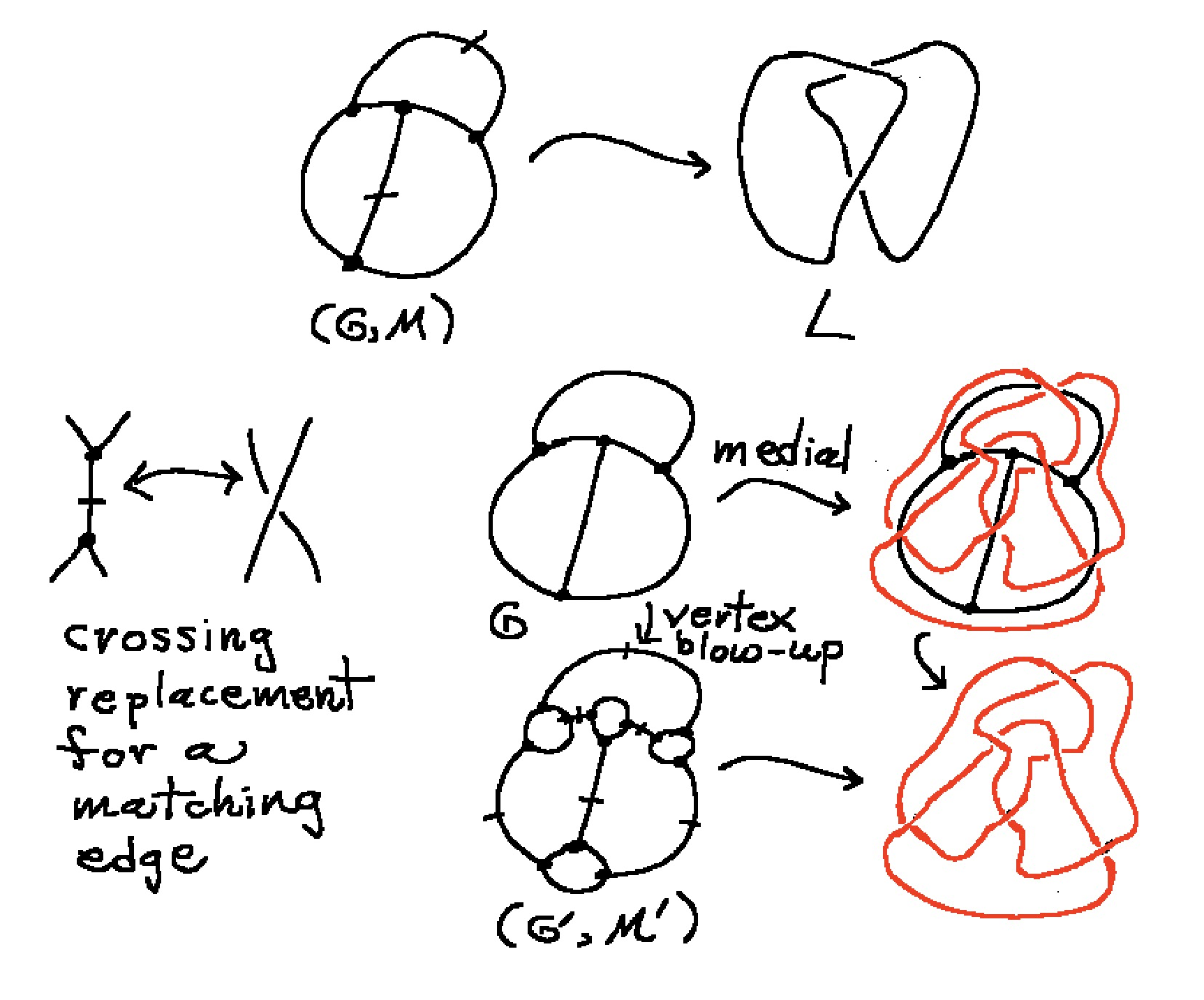}
    
     \caption{\bf Medial Graph and Blowing Up Nodes.}
     \label{medial}

\end{figure}

\begin{figure}
    
     \includegraphics[width=9cm]{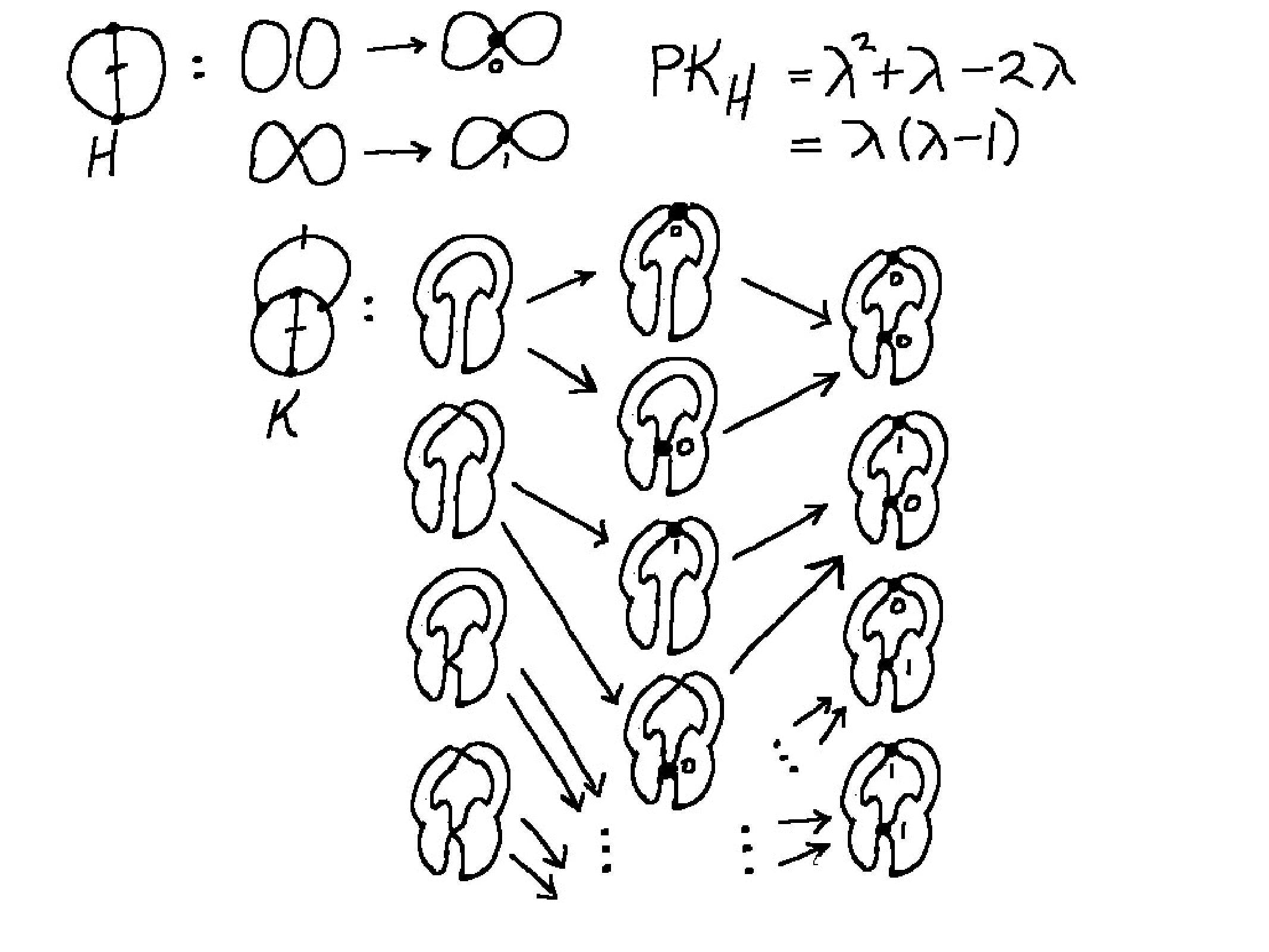}
     
     \caption{\bf Cubic Complex}
     \label{cubic}

\end{figure}

\begin{figure}
    
     \includegraphics[width=9cm]{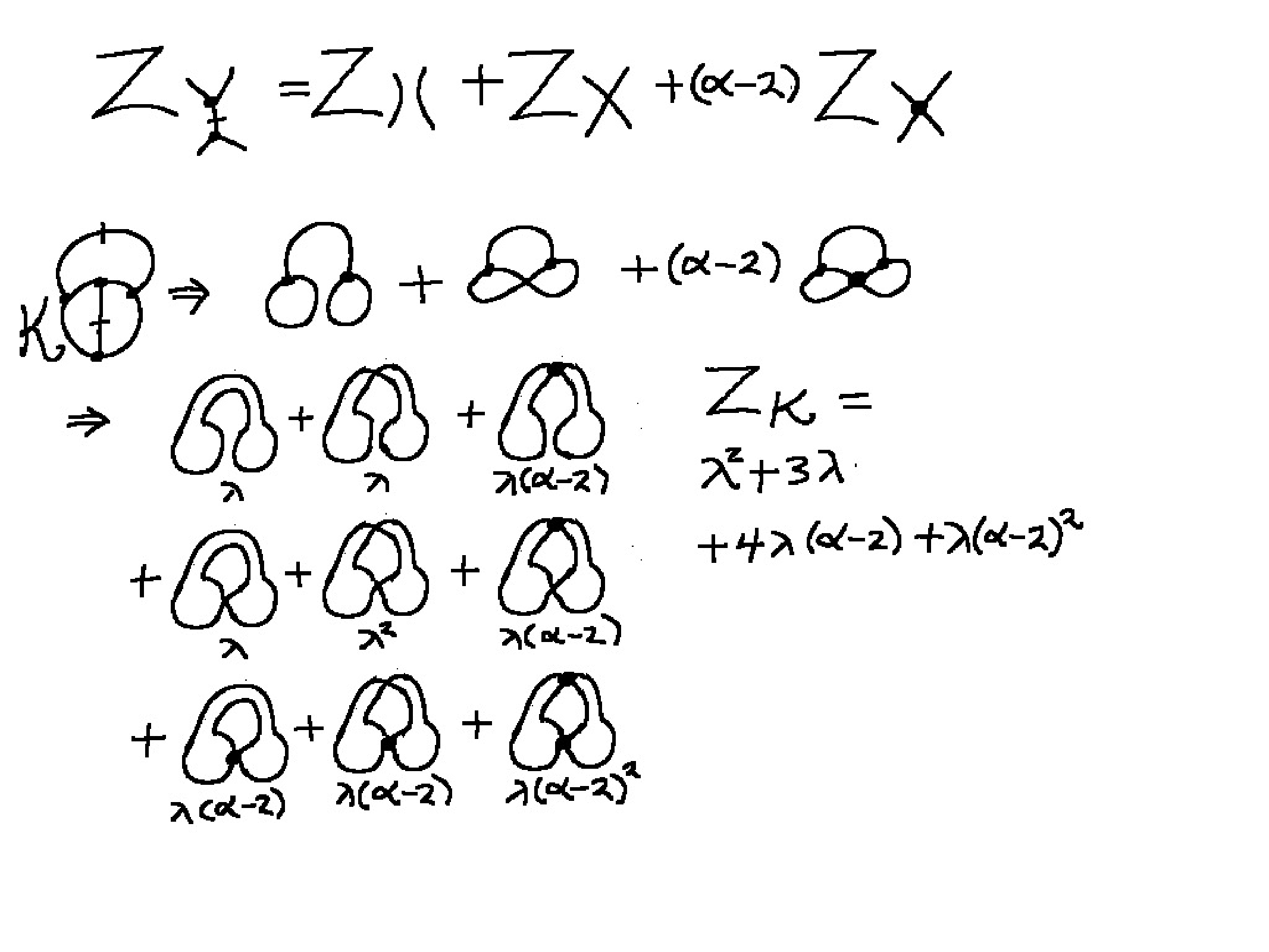}
    
     \caption{\bf Cubic Dichromatic}
     \label{cubicdichromatic}

\end{figure}

\begin{figure}[htb]
     
     \includegraphics[width=8cm]{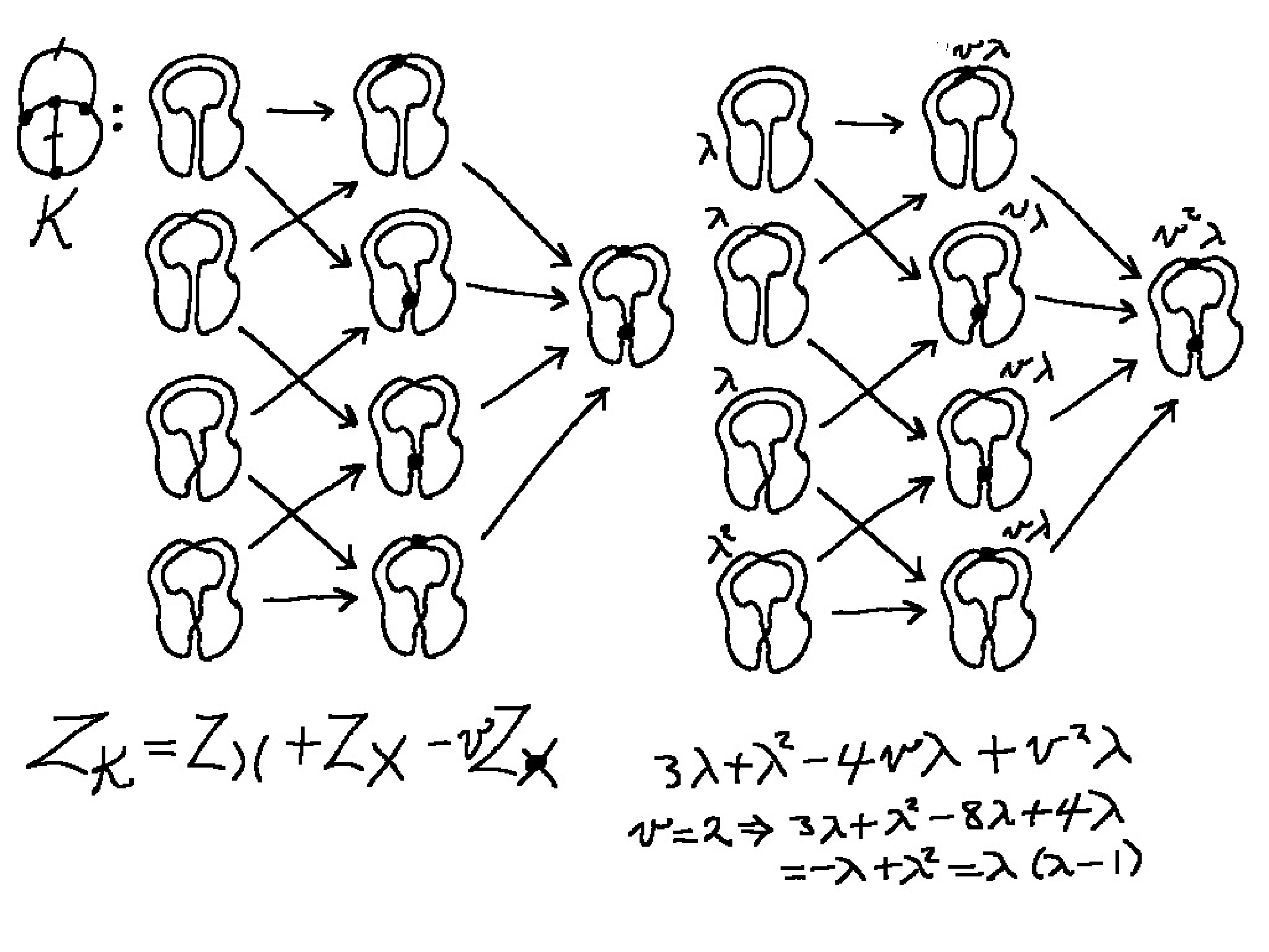}
    
     \caption{\bf Cubic Complex and Dichromatic Polynomial}
     \label{cubiccomplex}

\end{figure}

\section{Penrose Polynomials}
In this section we will consider the Penrose-Kauffman polynomial $PK[G,M]$ for any trivalent graph $G$ with a perfect matching $M.$  This polynomial is a generalization of the Penrose polynomial defined by 
Martin Aigner \cite{Aigner} (see also \cite{EM}) which is, in turn, a generalization of the Penrose evaluation $[G]$ for planar trivalent graphs \cite{P,BKM,LK4}. The original Penrose evaluation counts the number of {\it proper edge 3-colorings} of the graph $G$ where 
proper means that three distinct colors occur at each node of the graph, and for 3-colorings, only three colors are available.  In Figure~\ref{penrose} we indicate the generalization to $n$ colors that we shall use for a graph $G$ with a {\it perfect mathching} by which we mean a disjoint subset of edges of $G$ whose nodes comprise all the nodes of $G.$ The expansion formula shown in the figure is to be performed at matching edges.
The matching edges are designated by a transverse mark as indicated in the figure. The coloring convention is indicated as follows: Consider the glyph shown below and also in the figure.
$$ \YMTens$$
The glyph indicates four possible colors $a,b,c,d$ at the ends of the matching edge, with $a,b$ at one end and $c,d$ at the other end.
The {\it coloring condition} is that  $$a \ne b, c \ne d, and \,\, \{a,b\} = \{c,d\}.$$  Thus distinct colors occur at one end of the matching edge and the same two distinct colors occur at the other end.
This means that the two colors may or may not be permuted across the edge. This is the meaning of the two terms with parallel and crosed arcs in the expansion equation for the Penrose-Kauffman polynomial given below.
These terms are analogous to deletion in the Whitney formula discussed in the previous sections. The arcs can be freely colored either with the same color or with different colors. The crossed arcs with a dark bullet placed where they cross as in the glyph shown in the figure and in the formulas below are the analog of contraction in the Whitney formula. The dark bullet node means that {\it all four possible colors at the bullet node are equal to one another}. The analogy of the Whitney formula is as given below where we use $\lambda$ as the algebraic variable corresponding to $n$ colors. It is understood that the value of a disjoint loop is $\lambda$ and that the value of 
$k$ overlapping loops (transversely crossing each other without any nodes at the crossings) is $\lambda^{k}.$ We assume that $G$ is a cubic graph equipped with a perfect matching.
$$PK[ \YMGlyph ] = PK[ \CrossStarGlyph ] = PK[ \VDiag \ + PK[ \CDiag ] - 2PK[ \CDotDiag ],$$
$$PK[ G O ] = \lambda PK[G],$$
$$PK[ O ] = \lambda.$$\\
In Figure~\ref{penrose} we use $\lambda$ as an algebraic variable that can be specialized to a natural number $n.$ Note that the three-term formula above counts the number of colorings of the perfect matching graph because the parallel and crossed terms correspond to ``no local restriction on the coloring" (deletion) and the subtraction of twice the bullet node part of the formula  corresponds to subtracting the cases where colors are the same (contraction). Thus, by applying the formula, we count the cases where there are distinct colors in the sense of the conventions described above. The the factor of two in the formula is needed because a bullet node can come from either parallel arcs or crossed arcs.\\

In the Figure~\ref{penrose} we illustrate the coloring conditions for cubic graphs, and we illustrate how a matching edge can be replaced by a knot-theoretic crossing notation. The type of crossing determines, by the convention implicit in this figure, the expansion into parallel arcs - with the correct choice of parallels. Translating a perfect matching graph into a knot or link diagram gives a useful notation and raises some questions about the relationships of this approach to coloring polynomials with the topology of knots and links. Since the graph polynomials are not themselves invariant under isotopy moves for the topology, the relationship is not a direct one.
In \cite{LK4,KSW} the reader can find an exploration of properties of Penrose polynomials in the language of knot and link diagrams. In Figure~\ref{medial} we show the translation of a perfect matching graph $(G,M)$ into a link diagram. Note that topologically this link is isotopic to two disjoint circles. Nevertheless, the Penrose polynomial is non-trivial. In the lower part of the figure we show how the medial construction on a graph $G$ without perfect matching gives the same link diagram as the result of blowing up each trivalent vertex of the graph to a circle graph with three edges, then taking the edges of the original graph as matching edges for the new graph $G'.$ One sees directly that the link diagram for $G'$ is the same as the medial link construction for $G.$ In the medial construction on $G$ one makes a crossing for each edge of $G$ just as one makes a crossing for each matching edge of a perfect matching graph. Rather than give more detail about this construction, we refer the reader to what is shown in this figure.\\

In order to formulate a homology theory for the $PK$ polynomial, we consider all the graphs that are obtained by the analogues of deletion and contraction as described above. This means that we take the original graph and replace all perfect matching edges with either parallel or crossed arcs, or nodal crossed arcs. With appropriate enhancement labels, these will be the generators of the chain complex.  Call these configurations the {\it states}
$S$ of the graph $(G,M).$ The {\it loops} of a state $S$ without any nodal crossings consist in the closed paths in the state that are obtained by walking along the edges of the state and crossing transversely the arc at any crossing that one meets. For example in Figure~\ref{cubic} the left column of states have no crossing nodes and only the state at the bottom has two loops. If a state has crossing nodes, we say that two loops are {\it part of the same component} if an arc of one state makes a nodal crossing with an arc of the other state. (Note the analogy with contraction joining two components of a graph.) Thus we can count the number of components of a given state $S$ of $(G,M)$ by first counting the loops and then noting which loops are joined together via crossing nodes. For example, in the second column of states in Figure~\ref{cubic}  the last state, with one node, has one component because the two loops in this state are joined by that node. We let $c(S)$ denote the number of components of a state $S.$\\

The recursion formulas combined with the concept of states give the following formula for the $PK$ polynomial, where $n(S)$ is the number of nodal crossings in $S$ and $c(S)$ is the number of components of $S.$
$$PK_{G}(\lambda) = \sum_{S} (-1)^{n(S)} 2^{n(S)} \lambda^{c(S)}.$$
Note that this formula is analogous to the formula for the chromatic polynomial, where the states $S$ are subgraphs of the given graph. We will not make further connections here, but there is a natural graph associated to the states used in this section. Each loop in a state is associated with a vertex and each crossing node is an edge between the vertices corresponding to the loops at that node. For categorification, the factor of $2^{n(S)}$ can be taken to indicate that we take $2^{n(S)}$ copies of any given state with $n(S)$ nodes. See Figure~\ref{cubic} for an illustration where we label each crossing node with $0$ or $1$ according to whether it is regarded as coming from parallel arcs or from  a crossing. That is we associate an arrow from a state $S$ to a state $S'$ when $S'$ acquires a crossing node at one of its sites. Thus the arrow formally goes from parallel arcs to a crossing node labeled $0$, and from crossing arcs to a crossing node labeled $1.$  With the states so labeled, there will be $2^{k}$ states for each situation with $k$ crossing nodes. A partial drawing of the arrow diagrams for a given graph is shown in Figure~\ref{cubic}. In this form the expansion formula for the polynomial can be replaced with the labeled formula below
$$PK[ \YMGlyph ] = PK[ \CrossStarGlyph ] = PK[ \VDiag \ + PK[ \CDiag ] - PK[ \CDotDa ]-PK[ \CDotDb ].$$
When this formula is expanded into a state sum we sum over states with labels of $0$ or $1$ at each crossing node and can understand that there are arrows from states with parallel arcs to the $0$-crossing nodes and arrows from crossed arcs to the $1$-crossing nodes.\\

In order to categorify this formula we use enhanced states just as before. Thus we take $\lambda = (1 + q)$ and define an {\it enhanced state} $s$ corresponding to a state $S$ to be a labeling of the components of $S$ with the labels $x$ or $1$ with these labels regarded as elements of the ring $Z[x]/(x^2).$ Then we have (as before)  that 
$$PK_{G}(q) = \sum_{s} (-1)^{n(s)} q^{j(s)}$$
where the sum is over the enhanced states with $0,1$ labels at the nodal crossings, and $j(s)$ is the number of $x$ labels in the enhanced state $s.$ The chain complex ${ C}$ is generated by the enhanced states with ${ C}^{k,j}$ denoting the chains generated by enhanced states with $k$ crossing nodes and $j(s)$ labels of form $x.$  These labelings are in one-to-one correspondence with the terms in the expansion of $(1 + q)^{c(S)}.$ Signs are assigned in the differential, as before, by choosing an ordering of the matching edges and taking the corresponding ordering for the sites in each state. The differential $d: { C}^{k,j} \longrightarrow { C}^{k+1,j},$ preserves the count $j(s)$ just as in the case of the chromatic polynomial.
The differential is defined on enhanced states as a sum of signed partial differentials  $\partial(s) = s'$ for enhanced states $s$ and $s'$ so that $s'$ is obtained from $s$ by changing one parallel arc site or one crossing site to a nodal crossing site, labeling the result of joining two components by the product of their labels in $Z[x]/(x^2).$ The sign for a given partial differential in the summation for $d$ is equal to $(-1)^{\pi(s)}$ where $\pi(s)$ is the parity
of the number of sites preceding (in the given ordering) the site being operated on, that are not nodal crossing sites. It follows from this definition that $j$ is not changed by the differential and that $d^{2}=0.$\\

With this categorification we have 
$$PK_{G}(q) = \sum q^{j} \chi(H^{\star, j}(G)),$$
and all coefficients of the powers of $q$ are expressed as Euler characteristics in homology.\\

We formulate a dichromatic Penrose polynomial for $G$ with perfect matching $M$ by replacing the $2$ in the three term formula by a variables $v$ and $w$ for crossing nodes of type $0$ and $1$ respectively. Thus the equations for this dichromatic polynomial are as shown below.
$$Z[ \YMGlyph ] = Z[ \VDiag \ + Z[ \CDiag ] - vZ[ \CDotDa ] - wZ[ \CDotDb ],$$
$$Z[ G O ] = n Z[G],$$
$$Z[ O ] = n.$$\\
Note that $PK[G,M]$  and Z[G,M] are defined independent of any planar immersion of the graph $G.$  We have separated the variables $v$ and $w$ in order to articulate the state summation and for categorification.
The polynomial itself for a given graph, does not see the $0$ and $1$ labelings of these states and we have the equivalent formula
$$Z[ \YMGlyph ] = Z[ \VDiag \ + Z[ \CDiag ] - (v+w)Z[ \CDotDiag ],$$ for the graph polynomial. For computing the polynomial itself it is natural to replace $v+w$ by the single variable $v$ and proceed as illustrated in Figure~\ref{cubiccomplex} using $$Z[ \YMGlyph ] = Z[ \VDiag \ + Z[ \CDiag ] - vZ[ \CDotDiag ].$$ We can also use a special coefficient $\alpha$ and the expansion
$$Z[ \YMGlyph ] = Z[ \VDiag \ + Z[ \CDiag ] + (\alpha - 2)Z[ \CDotDiag ]$$ so that this version of the dichromatic Penrose polynomial is seen to have coefficients of powers of $\alpha$ that are impropriety
coloring polynomials analogous to the impropriety chromatic polynomials of the previous sections of this paper. Coloring at a matching edge will be called {\it improper} if all the colors that impinge on that edge are the same. This corresponds to equal colors on parallel or crossing arcs in the states. We will explore these variations of the Penrose dichromatic polynomial in a subsequent publication.\\

Using this formulation of the dichromatic Penrose polynomial for perfect matching graphs, we are in the position to repeat our definition of dichromatic homology in this context. 
Using the $0$ and $1$ labeled states $S$, we have $$Z[ \YMGlyph ] = Z[ \VDiag \ + Z[ \CDiag ] - vZ[ \CDotDa ] - wZ[ \CDotDb ]$$ and $$Z[G] = \sum_{S} (-1)^{n(s)} v^{n_{0}(S)}w^{n_{1}(S)}\lambda^{c(S)}$$ where $n(S) = n_{0}(S) + n_{1}(s)$ and $n_{0}(S)$ is the number of crossing nodes in $S$ labeled $0$ while $n_{1}(S)$ is the number of crossing nodes in $S$ labeled with $1.$ Letting $\lambda = 1 + q$, $v = 1 + p$ and $w = 1 + r,$  we can rewrite as a sum over enhanced states with $0,1$ labeled crossing nodes and new labels for the enhancements at crossing nodes. We choose variables $y$ and $1$ for the $0$ nodes and variables $z$ and $1$ for the $1$ nodes. Thus the expansion of $(1+p)^{n_{0}(s)}$ will correspond to the distribution of $1$ and $x$ on the $0$ crossing nodes, and the expansion of $(1+r)^{n_{1}(s)}$ will correspond to the distribution of $1$ and $z$ on the $1$ nodes. We use the ring $Z[x]/(x^2)$ as before for the component labeling of the enhanced states. Here we let $i(s)$ denote the number of components labeled $x$, $j(s)$ the number of crossing nodes of type $0$ labeled $y$ and $k(s)$ the number of crossing nodes of type $1$ labeled $z.$
Differentials are defined as before. Here $i$, $j$ and $k$ do not change under the differential and so we have a triply graded complex and a triply graded homology.
With this we have 
$$Z[G ] = \sum_{S} (-1)^{n(S)} (1+ q)^{c(S)}(1+p)^{n_{0}(S)}(1+r)^{n_{1}(S)}\lambda^{c(S)} $$
$$= \sum_{i,j,k} q^{i} p^{j}r^{k}\sum_{n} (-1)^{n} \sum_{s: n = n(s), i(s) = i, j(s) = j, k(s) = k} 1$$
$$= \sum_{i,jk} q^{i}p^{j}r^{k} \sum_{n} (-1)^{n} rank({ C}^{n,i,j,k}(G)$$
$$= \sum_{i,j,k} q^{i} p^{j}r^{k}\chi(H^{\star,i,j,k}(G))$$

The properties of this triply graded homology will be explored in subsequent work.\\

\section{Color Algebra - Another approach to categorification.}

In this section we return to the chromatic evaluation at a given integer $n$ (that is, there are $n$ colors)  so that 
$$C_{G}(n) = \sum_{S} (-1)^{e(S)}n^{c(S)}$$ where S runs over all subgraphs of $G$ (from a collection of disjoint nodes to the full graph $G$ with
$e(S)$ the number of edges in $S$ and $c(S)$ the number of connected components in $S.$ This formula would be the Euler characteristic of a complex ${ C}(G)$
if the complex is graded according to the number of edges in the subgraph so that ${ C}^{k}(G)$ denotes the module generated by subgraphs with $k$ edges, with each component labeled with one color from the set of possible colors. See Figure~\ref{graphcomp} for an illustration of some graph complexes and the alternating rank computation that gives both the Euler characteristic and the color count.\\

With these assumptions we would have that  $n^{c(S)}$ is the number of module generators corresponding to the state $S$ and
$$C_{G}(n) = \sum_{S} (-1)^{e(S)}n^{c(S)} = \sum_{k} (-1)^{k} \sum_{S: c(S)=k} n^{k} =\sum_{k} (-1)^{k} rank({ C}^{k}(G)) = \chi({ C}(G)).$$
Given that there is a homology associated with the complex, we then have $C_{G}(n)=  \chi({ H}(G)).$
In this form of categorification, certain numbers appear as Euler characteristics, and we have a categorification of the chromatic evaluation itself.
\\

In order to realize the categorification described above, we choose an algebra structure on the $n$ colors.
We take commuting algebraic variables $\{ x_1, x_2, \cdots x_n \}$ as generators of an associative algebra over the integers. Any associative algebra will suffice just so long as $x_{i}x_{j}$ is equal to some $x_{k}$ or $0$ for each choice of $i$ and $j.$ For example we can take  $x_{i} x_{j} = 0$ if $i \ne j$ and $x_{i}^{2}= x_{i}.$ Call this algebra ${ A}(n).$ We shall refer to the $x_i$ as {\it colors}. If $S$ is a subgraph of $G$ let $s$ denote an {\it enhancement} of $S$ obtained by assigning a color to each component of $S.$ Define partial boundaries $\partial_{e}(s)$ for each enhanced state $s$ and an edge $e$ in $G - S$ so that 
\begin{enumerate}
\item Suppose the endpoints of $e$ are in the same component of $S$ and this component is colored $x.$ Then $\partial_{e}(s) = s'$ where $s'$ is obtained from $s$ by adding the one edge $e$ and taking the product color $x$ for the component containing $e.$ The other colors in $s'$ are inherited from $s.$ 
\item Suppose that the endpoints of $e$ are in different components of $s$ with colors $x$ and $x'.$ Then $\partial_{e}(s) = s'$ where $s'$ is obtained from $s$ by adding one edge $e$ and taking the product color $xx'$
 for the component containing $e.$ The other colors in $s'$ are inherited from $s.$ 
 \item  If the color product is $0$ then the partial boundary is $0.$
 \item Note that if $s$ has $k$ edges, then 
$\partial_{e}(s)=s'$ has $k+1$ edges. Thus $$\partial_{e}:{ C}^{k}(G) \longrightarrow { C}^{k+1}(G)$$ by extending each such map linearly over the module generators. 
\item Choose an ordering of the edges of $G.$ Define $$d:{ C}^{k}(G) \longrightarrow { C}^{k+1}(G)$$ to be the signed sum of $\partial_{e}$ over all $e$ in $G-S$ 
$$d(s) = \sum_{e \in G-S} (-1)^{\pi(e,s)} \partial_{e}(s)$$ where $\pi(e,s)$ denotes the number of edges in $G$ preceding $e$ in the chosen ordering that are in $G-S$ and $e \in G-S$ means that $e$ is an edge in the complement of $S.$ It then follows that $d \circ d = 0,$ so that there is a chain complex and well defined homology.
\end{enumerate}

Using the algebra ${ A}(n),$ let $H^{k}(G, { A}(n))$ denote the homology defined above at $k$ edges. Thus we have $$C_{G}(n) = \chi(H^{\star}(G, { A}(n))).$$
Let $G$ have vertex set $\{ v_1, v_2, \cdots , v_m\}.$
Suppose that  $s$ is any assignment of colors $s(i) = s(v_i)$  to the nodes $v_i$ with values in the set of algebra generators $\{x_1, x_2, \cdots , x_n \}.$ By our definitions, $s$ is an enhanced state in $C^{0}(G, { A}(n)).$
Since $x_i x_j = 0$ when $i \ne j,$ we see that for the differential $d$ in the chain complex that we have defined, $d(s) = 0$ if and only if $s(i) \ne s(j)$ whenever $v_i$ and $v_j$ are distinct vertices that are endpoints of an edge in $G$. In this case the state
$s$ is a proper coloring of $G.$ It follows from these remarks that \\

\noindent{\bf Proposition.} Using the algebra ${ A}(n),$ we have that $H^{0}(G, { A}(n))$ is freely generated by the proper colorings of the graph $G.$ Hence $ rank(H^{0}(G, { A}(n))) = C_{G}(n).$
The higher indexed homology groups  $H^{\star}(G, { A}(n))$ are zero.\\

With these definitions we have a method to categorify chromatic evaluations. There are many choices of algebra ${A}.$ For example, let $H$ be any finite abelian (we did not specify order of multiplication in the boundary maps) group of order $n$ and let ${ A} = Z[H]$ be the integral group ring of $H.$ This is an associative algebra, and we take the elements of $H$ as the colors. Then for any graph $G$ there is a homology theory $H^{\star}(G, Z[H])$ with $C_{G}(n) = \chi(H^{\star}(G, Z[H])).$\\

If $n=4$ then it is natural to take the elements $\{ 1,a,b,c \}$ of the Klein $4$-group as the set of colors and let ${ K}(4)$ be the group ring over the integers of the Klein 4-group.  In the Klein $4$-group we have that $1$ is the identity element, $a^2 = b^2 = c^2 = 1$ and $ab=ba= c, ac = ca = b, bc = cb = a.$ With this color algebra in mind, we have, by the four color theorem, that for any $1$-connected plane graph $G,$  the chromatic evaluation $C_{G}(4)$ is non-zero and hence $\chi(H^{\star}(G,{ K}(4))) \ne 0.$ Thus the homology $H^{\star}(G,{ K}(4) )$ is non-trivial. If there were an alternate route to proving that this homology is sufficiently non-trivial, then the four color theorem would follow.\\

See Figure~\ref{graphcomp}. For the graph complexes in this figure it is not hard to see that the color count using the Klein $4$-group is concentrated at the $H^{0}(G).$ For more complex graphs it remains to be seen if the Klein $4$-group homology will be distributed in different dimensions in the chain complex.\\

\noindent{\bf Remark.} In this paper we have defined a number of homology groups related to graphs and have indicated how these categorify chromatic polynomials, dichromatic polynomials and Penrose polynomials.
In this final section we have pointed out that these ideas also lead categorifications of numerical chromatic evaluations, identifying them as Euler characteristics of a panoply of homology theories. All of these situations will be explored further structurally and with specific computations.\\

\end{document}